\documentclass[12pt,reqno]{amsart}
\usepackage[utf8x]{inputenc}

\usepackage{hyperref}
\hypersetup{
    colorlinks=true,       
    linkcolor=blue,          
    citecolor=blue,        
    filecolor=blue,      
    urlcolor=cyan  
}
\usepackage{latexsym,amssymb,amsmath,amsthm,bbm}
\usepackage{epsfig,mathrsfs}
\usepackage{pst-eucl}
\usepackage{pst-plot,pstricks-add}

\newcommand{\R}{{\ensuremath{\mathbb{R}}}}
\newcommand{\N}{{\ensuremath{\mathbb{N}}}}

\renewcommand{\P}{\ensuremath{\mathbb{P}}}
\renewcommand{\dj}{d\kern-0.4em\char"16\kern-0.1em}
\newcommand{\E}{\ensuremath{\mathbb{E}}}

\newcommand{\sig}{\ensuremath{\mathcal{F}}}

\newtheorem{Thm}{Theorem}[section]

\newtheorem{Lem}[Thm]{Lemma}
\newtheorem{Prop}[Thm]{Proposition}

\theoremstyle{remark}
\newtheorem{Rem}[Thm]{Remark}

\theoremstyle{definition}

\theoremstyle{definition}

\theoremstyle{definition}

\begin{document}
\numberwithin{equation}{section}
\bibliographystyle{amsalpha}

\title[Harmonic functions of L\' evy processes]{On harmonic functions of
symmetric L\' evy processes}

\author{Ante Mimica}
\address{Fakult\"{a}t f\"{u}r Mathematik\\Universit\"{a}t Bielefeld\\Germany}
\email{amimica@math.uni-bielefeld.de}
\thanks{ Research supported in part by German Science Foundation DFG via IGK
"Stochastics and real world models" and SFB 701.}
\date{}

\maketitle

\begin{abstract}
We consider some classes of L\' evy processes for which the
estimate of Krylov and Safonov (as in \cite{BL1}) fails and thus it is not
possible to use the standard iteration technique to obtain a-priori H\"
older continuity estimates of harmonic functions. Despite the faliure of this
method, we obtain some a-priori
regularity estimates of harmonic functions for these processes. Moreover, we
extend results
from \cite{SSV2} and obtain asymptotic behavior of the Green function and the
L\' evy
density for a large class of subordinate Brownian motions, where the Laplace
exponent of the corresponding subordinator is a slowly varying function. 
\end{abstract}

\section{Introduction}

Recently there has been much interest in investigation of the continuity
properties of
harmonic functions with the respect to various non-local operators. 
An example of such operator $\mathcal{L}$ is of the form
\begin{equation}\label{eq:tmp-l-op}
	(\mathcal{L}f)(x)=\int_{\R^d\setminus\{0\}}\left(f(x+h)-f(x)-\langle
\nabla f(x),h\rangle \mathbbm{1}_{|h|\leq 1}\right)n(x,h)\,dh
\end{equation}
for $f\in C(\R^d)$ bounded. Here $n\colon
\R^d\times(\R^d\setminus\{0\})\rightarrow [0,\infty)$ is a measurable function
satisfying
\[
	c_1|h|^{-d-\alpha}\leq n(x,h)\leq c_2|h|^{-d-\alpha}\,,
\]
for some constants $c_1,c_2>0$ and $\alpha\in [0,2]$. 

It is known that for $\alpha\in (0,2)$ H\" older regularity estimates hold for $\mathcal{L}$-harmonic functions (see \cite{BL1} for a
probabilistic and
\cite{Si} for an analytic approach) .

For $\alpha=0$, techniques developed so far are  not applicable. One of our
aims is to investigate this case using a probabilistic approach. 

In many cases the operator of the form (\ref{eq:tmp-l-op}) can be understood as the infinitesimal generator of a Markov jump process. The kernel $n(x,h)$ can be thought of as the measure of intensity of jumps of the process. 

Let us describe the stochastic process we are considering. Let $S=(S_t:t\geq
0)$ be a subordinator such that its Laplace exponent $\phi$ defined by $\phi(\lambda)=-\log\left(\E e^{-\lambda S_1}\right)$ satisfies
\begin{equation}\label{eq:intro-sw-13}
			\lim_{\lambda\to+\infty }\frac{\phi'(\lambda
x)}{\phi'(\lambda)}=x^{\frac{\alpha}{2}-1}\ \ \text{ for any }\ \ x>0
\end{equation}
for some $\alpha\in [0,2]$ . Let $B=(B_t,\P_x)$ be an independent Brownian
motion in $\R^d$ and define a new process $X=(X_t,\P_x)$ in $\R^d$ by $X_t=B(S_t)$ . It 
is called the subordinate Brownian motion. 

{\bf Example 1:} Let $S$ be a subordinator with the Laplace exponent $\phi$
satisfying
 \[\lim_{\lambda \to +\infty }\frac{\phi(\lambda)}{
\lambda^{\alpha/2}\ell(\lambda)}=1\] with
$\alpha\in (0,2)$ and $\ell\colon (0,\infty)\rightarrow (0,\infty)$ that varies
slowly at infinity (i.e. for any $x>0$,  $\frac{\ell(\lambda
x)}{\ell(\lambda)}\to 1$ as $\lambda\to+\infty$).\\
Then (\ref{eq:intro-sw-13})
holds (we just use Theorem 1.7.2 in
\cite{BGT} together with the fact that  $\phi'$ is
decreasing and $\phi(\lambda)=\int_0^\lambda \phi'(t)\,dt$) and 
\[
n(x,h)=n(h)\asymp |h|^{-d-\alpha}\ell(|h|^{-2}),\ |h|\to 0+\,,
\]
which means that $\frac{n(h)}{|h|^{-d-\alpha}\ell(|h|^{-2})}$ stays between two
positive constants as $|h|\to 0+\,$.

Choosing  $\ell\equiv 1$ we see  that the rotationally invariant
$\alpha$-stable process (whose infinitesimal generator is the fractional
Laplacian $\mathcal{L}=-(-\Delta)^{\frac{\alpha}{2}}$) is  included in this
class.  Other choices of $\ell$ allow us
to consider processes which are not invariant under time-space scaling.

The processes described in Example 1 with $\ell\equiv1$ belong to the
class of L\' evy stable or, more generally, stable-like Markov jump processes.
The potential
theory of such processes is well
investigated (see \cite{BL1,CK2,SV1,BS,BKa,RSV,KS,CK,Mi,Sz2,KM}). 
For example, it is known that harmonic functions of
 such processes satisfy H\"{o}lder regularity estimates and the scale invariant
Harnack inequality holds. Also, two sided heat kernel estimates are obtained for
these processes. 

Not much  is known about harmonic functions in the case when
the corresponding subordinators belong to 'boundary' cases, i.e. $\alpha\in
\{0,2\}$ (see \cite{SSV2,Mi2,Mi3}). For
the class of geometric stable processes only the non-scale invariant Harnack
inequality was proved and on-diagonal heat kernel upper estimate is not finite
(see \cite{SSV2}).

The following two examples belong to these 'boundary' cases and are covered by our approach.  

{\bf Example 2:} Let $\phi(\lambda)=\log(1+\lambda)$. The corresponding  process $X$ is known as the variance
gamma process. In (\ref{eq:intro-sw-13}) we have $\alpha=0$ . Moreover, it
will be proved (see Theorem \ref{thm:jump_asymp}) that
\begin{equation}\label{eq:tmp_geom_st_1001}
	n(x,h)=n(h)\asymp |h|^{-d},\ |h|\to 0+\,.
\end{equation}
This example can be generalized in various ways. For example, we can take $k\in
\N$ and consider $\phi_k=\underbrace{\phi\circ\ldots\circ\phi}_{k\
\text{times}}$. Then (see Theorem \ref{thm:jump_asymp} or Subsection
\ref{subsec:geom_stable})
\[
	n_k(x,h)=n_k(h)\asymp |h|^{-d}\left(\underbrace{\log\cdots\log}_{k-1
\text{ times }} \tfrac{1}{|h|}\cdot\ldots\cdot  \log\log \tfrac{1}{|h|}\cdot
\log\tfrac{1}{|h|}\right)^{-1},\ |h|\to 0+\,.
\]
Another generalization is to consider $\phi(\lambda)=\log(1+\lambda^{\beta/2})$
for some $\beta\in (0,2]$. The process $X$ is known as the $\beta$-geometric
stable process and the behavior of $n$ is given also by
(\ref{eq:tmp_geom_st_1001}).

{\bf Example 3:} Let $\phi(\lambda)=\frac{\lambda}{\log(1+\sqrt{\lambda})}$ .
Then
in (\ref{eq:intro-sw-13}) we have  $\alpha=2$ and  (see Theorem \ref{thm:jump_asymp})
\[
	n(x,h)=n(h)\asymp |h|^{-d-2}\left(\log\tfrac{1}{|h|}\right)^{-2}, \
|h|\to 0+\,.
\]
This behavior shows that small jumps of this process have higher intensity than
small jumps of any stable process. 

A measurable bounded function $f\colon \R^d\rightarrow \R$ is said to be harmonic in an
open set $D\subset \R^d$ if for any relatively compact open set $B\subset
\overline{B}\subset D$ 
\[
	f(x)=\E_x f(X_{\tau_B})\ \text{ for any }\ x\in B,
\]
where $\tau_B=\inf\{t>0\colon X_t\not\in B\}$ .


The main theorem is the following regularity result, which covers cases $
\alpha\in [0,1)$.  The novelty of this result rests on the case $\alpha=0$. By
$B_r(x_0)$ we denote the open ball with center $x_0\in \R^d$ and radius $r>0$ .

\begin{Thm}\label{tm:main}
	Let $S$ be a subordinator  such that its L\' evy and potential
measures have decreasing densities. Assume that the Laplace exponent of $S$ satisfies
(\ref{eq:intro-sw-13}) with $\alpha\in [0,1)$. Let $X$ be the corresponding
subordinate Brownian motion and let $d\geq 3$. 

There is a constant $c>0$ such
that
for any
\mbox{$r\in (0,\frac{1}{4})$} and any bounded function $f\colon \R^d\rightarrow
\R$ which is
harmonic in $B_{4r}(0)$,
	\[
		|f(x)-f(y)|\leq c\|f\|_\infty
\frac{\phi\left(r^{-2}\right)}{\phi\left(|x-y|^{-2}\right)}\ \ \text{ for all }\
x,y\in B_{\frac{r}{4}}(0)\,.
	\]
\end{Thm}

Applying Theorem \ref{tm:main} to Example 1, we obtain expected H\" older regularity
estimates. Within  this example the result is new when the scaling is lost, e. g.
\[\phi(\lambda)=\lambda^{\frac{\alpha}{2}}\left[\log(1+\lambda)\right]^{
1-\frac{\alpha}{2}}\,.\]

The situation is more interesting in Example 2, e.g. for  the geometric stable
process. For this process we obtain logarithmic
regularity estimates:
\[
		|f(x)-f(y)|\leq c\|f\|_\infty
\log(r^{-1})\frac{1}{\log(|x-y|^{-1})}\,,
\]
It is still unknown whether H\"{o}lder regularity estimates hold  for harmonic
functions of this process (or, generally, of the
processes belonging to the case $\alpha=0$).

Let us explain why known analytic and probabilistic techniques do not work in
the case $\alpha=0$. 
The main idea in the proof of the a priori H\" older estimates of harmonic
functions relies on the estimate of Krylov and Safonov. 

In probabilistic setting this estimate can be formulated as follows. 
There is a constant $c>0$ such that for every closed subset $A\subset B_r(0)$
and $x\in B_{\frac{r}{2}}(0)$
\begin{equation}\label{eq:kr_saf}
	\P_x(T_A<\tau_{B(0,r)})\geq c\,\frac{|A|}{|B_r(0)|}\,,
\end{equation}
where $T_A=\tau_{A^c}$ is the first hitting time of $A$  and $|A|$ denotes the Lebesgue measure of the set $A$ .

Performing a computation similar to the one in the proof of Proposition 3.4
in 
\cite{BL1} (see also Lemma 3.4 in \cite{SV1}) we deduce
\[
	\P_x(T_{A}<\tau_{B_r(0)})\geq
c\,\frac{r^{-2}\phi'(r^{-2})}{\phi(r^{-2})}\frac{|A|}{|B_r(0)|}\, .
\]

If $\alpha\in (0,2)$, it can be seen that
$\frac{r^{-2}\phi'(r^{-2})}{\phi(r^{-2})}\asymp 1$ as $r\to 0+$ . This gives estimate of the form (\ref{eq:kr_saf}) and thus the
standard Moser's iteration procedure  for obtaining a-priori H\" older
regularity estimates of
harmonic functions can be applied (see the proof of Theorem 4.1 in \cite{BL1}
for a probabilistic version).

The situation is quite different for $\alpha=0$. To find a counterexample we
will use the following result. 

\begin{Prop}\label{prop:count-ks}
Let $S$ be as subordinator  such that its L\' evy and potential
measures have decreasing densities and whose Laplace exponent satisfies
(\ref{eq:intro-sw-13}) with $\alpha\in [0,1)$. Let $X$ be the corresponding
subordinate Brownian motion and let $d\geq 3$.

 There is a constant $c>0$ such
that
for every $r\in
(0,1)$ and $x\in B_{\frac{r}{4}}(0)$
\[
 \P_x(X_{\tau_{B_{\frac{r}{2}}(0)}}\in B_r(0)\setminus B_{\frac{r}{2}}(0))\leq
c\frac{r^{-2}\phi'(r^{-2})}{\phi(r^{-2})}\,.
\]
\end{Prop}

For $\alpha=0$ it will follow that  $\lim\limits_{r\to
0+}\frac{r^{-2}\phi'(r^{-2})}{\phi(r^{-2})}= 0$ (see (\ref{eq:de_haan_limit})).
Therefore in this case (\ref{eq:kr_saf}) does not hold, since
\[
 \lim_{r\to 0+}\P_0(T_{B_r(0)\setminus B_{\frac{r}{4}}(0)}<\tau_{B_r(0)})\leq
\lim_{r\to 0+}\P_0(X_{\tau_{B_{\frac{r}{4}}(0)}}\in
B_r(0)\setminus B_{\frac{r}{4}}(0))=0\,.
\]

Considering process $X$ in the setting of metric measure spaces (as in \cite{CK}
or \cite{BGK}) the news feature appears. Theorem \ref{thm:jump_asymp} shows
that the jumping kernel of the process $X$ is of the form $n(x,h)=
j(|h|)$ with
\[
	j(r)\asymp  \frac{r^{-2}\phi'(r^{-2})}{\phi(r^{-2})}\cdot\frac{\E_0
\tau_{B_r(0)}}{|B_r(0)|},\ r\to 0+\,.
\]
In the case $\alpha=0$ the term $\frac{r^{-2}\phi'(r^{-2})}{\phi(r^{-2})}$
becomes significant. This has not yet been treated within this framework.

The latter discussion shows that the question of the continuity of harmonic
functions becomes interesting even in the case when the kernel $n(x,h)$ is space
homogeneous, or in other words, in the case of a L\' evy process. There is no
known technique that covers this situation in the case of a more general jump
process.

Our technique relies on asymptotic properties of the underlying subordinator.
The potential density can be analyzed using the de Haan theory of slow variation
(see \cite{BGT}). 

On the other hand, there is no known Tauberian theorem that can be applied to
obtain asymptotic
behavior of the L\' evy density $\mu$  of the subordinator. For this purpose we
perform asymptotic inversion of the Laplace transform (see Proposition
\ref{prop:sub-mu1}) to get  
\[
\mu(t)\asymp t^{-2}\phi'(t^{-2}),\ t\to 0+\,.
\]

These techniques allow us to extend results from
\cite{SSV2} to much wider class of subordinators whose Laplace exponents are
logarithmic or, more generally, slowly varying functions. 

Although we do not obtain regularity estimates of harmonic functions for cases when
 $\alpha\in [1,2]$, it is possible to say  something about the behavior of the jumping
kernel and the Green function. In this sense, the case $\alpha=2$ is also new.
For example the Green function of the process corresponding to the Example 3
above has the following behavior:
\[
 G(x,y)\asymp |x-y|^{2-d}\log(|x-y|^{-1}),\ |x-y|\to 0+\,.
\]
We may say that such process $X$ is  'between' any stable process and
Brownian motion. 

Let us briefly comment the technique we are using to prove the regularity result.
In Section \ref{sec:prelim} it will be seen  that any bounded function $f$
which is harmonic in $B_{2r}(0)$ can be represented as
\[
 f(x)=\int_{\overline{B_r(0)}^c}K_{B_r(0)}(x,z)f(z)\,dz,\ x\in B_r(0)\,,
\]
where $K_{B_r(0)}(x,z)$ is the Poisson kernel of the ball $B_r(0)$. 

The following estimate of differences of Poisson kernel is the key to the proof of Theorem \ref{tm:main}: 
\[
 |K_{B_r(0)}(x_1,z)-K_{B_r(0)}(x_2,z)|\leq\left\{\begin{array}{cl}
c\,|z|^{-d}\frac{\phi((|z|-r)^{-2})}{\phi(|x-y|^{-2})} & r<|z|\leq 2r\\
\frac{j(\frac{|z|}{2})}{\phi(|x-y|^{-2})} & |z|> 2r\end{array}\right.
\]
for $x_1,x_2\in B_{\frac{r}{8}}(0)$ (see Proposition \ref{prop:poisson_k_est}).

Similar type estimate has been obtained in \cite{Sz2} for stable L\' evy processes
using scaling argument and the explicit behavior of the transition density. In our setting there are
many cases where the behavior of the transition density is not known and the scaling argument does not work. Our idea is to establish the following Green function difference estimates:
\[
 |G(x_1,y)-G(x_2,y)|\leq c\frac{r^{-2}\phi'(r^{-2})}{\phi(r^{-2})^2}\,
r^d\,\left(1\wedge \tfrac{|x-y|}{r}\right)
\]
for all $y\in \R^d$ and $x_1,x_2\not\in B_r(y)$ (see
Proposition \ref{prop:green_diff_1}).

The paper is organized as follows. In Section \ref{sec:prelim} we introduce all
concepts we need throughout the paper. Section \ref{sec:asymp} is devoted to the
study of subordinators. We obtain asymptotic properties of L\' evy and potential
densities. In Section \ref{sect:sbm} asymptotical properties of the  Green function and  L\' evy density of the subordinate Brownian motions are obtained.
Difference estimates of the Green function and the Poisson kernel are the main subject of Section
\ref{sect:diff}. This type of estimates are the main ingredient in the proof of
the regularity result in Section \ref{sect:main}. In Section
\ref{sec:examples} we
apply our results to some new examples.

{\bf Notation. }For two functions $f$ and $g$  we write $f\sim g$ if $f/g$
converges to 1 and $f\asymp g$ if $f/g$ stays between two positive constants.
The $n$-th derivative of $f$ (if exists) is denoted by $f^{(n)}$.

The logarithm with base $e$ is denoted by $\log$ and we introduce the following
notation for iterated logarithms: $\log_1=\log$ and $\log_{k+1}=\log\circ
\log_k$ for $k\in\N$.

We say that $f\colon\R\rightarrow \R$ is increasing if
$s\leq t$ implies $f(s)\leq f(t)$ and analogously for a decreasing function. 

The standard Euclidian norm and   the standard inner product
in $\R^d$ are
denoted by $|\cdot|$ and $\langle\cdot,\cdot\rangle$, respectively. By
$B_r(x)=\{y\in\R^d\colon |y-x|<r$\} we denote the open ball centered at $x$ with
radius $r>0$. The Gamma function is defined by $\Gamma(\rho)=\int_0^\infty
t^{\rho-1}e^{-t}\,dt$ for  $\rho>0$.

\section{Preliminaries}\label{sec:prelim}

\subsection{L\' evy processes and their potential theory}

A stochastic process $X=(X_t\colon t\geq 0)$ with values in  $\R^d$ $(d\geq 1)$ defined on a probability space
$(\Omega,\sig,\P)$  is said to be a L\' evy process if it
has independent and stationary increments, its trajectories are $\P$-a.s. right
continuous with left limits and $\P(X_0=0)=1$\,.

The characteristic function of $X_t$ is always of the form
\[
	\E \exp{\{i\langle \xi,X_t\rangle \}}=\exp{\{-t\Phi(\xi)\}},
\]
where $\Phi$ is called the characteristic (or L\' evy) exponent of $X$. It has
the following L\' evy-Khintchine representation
\[
	\Phi(\xi)=i\langle \gamma,\xi\rangle +\frac{1}{2}\langle A\xi,\xi\rangle
+\int_{\R^d}\left(1-e^{i\langle x,\xi\rangle}+i\langle
x,\xi\rangle\mathbbm{1}_{\{|x|\leq 1\}}\right)\Pi(dx).
\]
Here $\gamma\in \R^d$, $A$ is a non-negative definite symmetric $d\times d$ real
matrix and $\Pi$ is a measure on $\R^d$, called the L\' evy measure of $X$,
satisfying
\[
	\Pi(\{0\})=0\ \text{ and }\ \int_{\R^d}(1\wedge |x|^2)\Pi(dx)<\infty\,.
\]

The Brownian motion $B=(B_t\colon t\geq 0)$ in $\R^d$  with transition density
$p_0(t,x,y)=(4\pi t)^{-d/2}\exp{\left\{-\frac{|x-y|^2}{4t}\right\}}$ is an
example of a L\' evy process with the characteristic exponent
$\Phi(\xi)=|\xi|^2$ . 

A subordinator is a stochastic process $S=(S_t\colon t\geq 0)$ which is a L\'
evy process in $\R$ such that $S_t\in [0,\infty)$ for every $t\geq 0$. 
In this case it is more convenient to consider the Laplace transform of $S_t$:
\[
	\E \exp{\{-\lambda S_t\}}=\exp{\{-t\phi(\lambda)\}},\ \lambda >0.
\]
The function $\phi\colon (0,\infty)\rightarrow (0,\infty)$ is called the Laplace
exponent of $S$ and it has the following representation
\[
	\phi(\lambda)=\gamma \lambda + \int_{(0,\infty)}(1-e^{-\lambda
t})\mu(dt).
\]
Here $\gamma\geq 0$ and $\mu$ is also called the L\' evy measure of $S$ and it
satisfies the following integrability condition: $\int_{(0,\infty)}(1\wedge
t)\mu(dt)<\infty$.

The potential measure of the subordinator $S$ is defined by 
\[
U(A)=\E\left[\int_0^\infty \mathbbm{1}_{\{S_t\in A\}}\,dt\right]\ \text{
for a measurable }\ A\subset
[0,\infty)\,.	
\]
The Laplace transform of $U$ is then
\begin{equation}\label{eq:prelim-tmp1}
	\mathcal{L}U(\lambda):=\int_{(0,\infty)}e^{-\lambda
t}\,U(dt)=\frac{1}{\phi(\lambda)}\,.
\end{equation}

Assume that the processes $B$ and $S$ just described are independent. We define
a new stochastic process $X=(X_t\colon t\geq 0)$ by $X_t=B(S_t)$ and call it the
subordinate Brownian motion. It is a L\' evy process with the characteristic
exponent $\Phi(\xi)=\phi(|\xi|^2)$ and the L\' evy measure of the form
$\Pi(dx)=j(|x|)\,dx$ with
\begin{equation}\label{eq:prelim-jot}
	j(r)=\int_{(0,\infty)}(4\pi
t)^{-d/2}\exp{\left\{-\frac{r^2}{4t}\right\}}\,\mu(dt)\,.
\end{equation}

The process $X$ has the transition density and it is given by
\begin{equation}\label{eq:sbm_trans}
p(t,x,y)=\int_{[0,\infty)}(4\pi
s)^{-d/2}\exp{\left\{-\frac{|x-y|^2}{4s}\right\}}\P(S_t\in ds)\,.
\end{equation}

When $X$ is transient, we can define the Green function of $X$ by
\[
	G(x,y)=\int_{(0,\infty)}p(t,x,y)\,dt,\ x,y\in \R^d,\ x\not =y\,.
\]	
The Green function can be considered as the density of the Green measure defined
by
\[
	G(x,A)=\E_x\left[\int_0^\infty \mathbbm{1}_{\{X_t\in A\}}\,dt\right],\
A\subset \R^d \ \text{ measurable },
\]
since $G(x,A)=\int_A G(x,y)\,dy$.

Using (\ref{eq:sbm_trans}) we can rewrite it as $G(x,y)=g(|y-x|)$ with
\begin{equation}\label{eq:prelim-green}
	g(r)=\int_{(0,\infty)}(4\pi
t)^{-d/2}\exp{\left\{-\frac{r^2}{4t}\right\}}U(dt)\,.
\end{equation}

Let $D\subset \R^d$ be a bounded open set. We define the process killed upon exiting $D$ 
$X^D=(X_t^D\colon t\geq 0)$  by
\[
	X_t^D=\left\{\begin{array}{cl}
		X_t & t<\tau_D\\
		\partial & t\geq \tau_D,
	\end{array}\right.
\] 
where $\partial$ is an extra point adjoined to $D$.

Using the strong Markov property we can see that the Green measure of $X^D$ is
\begin{align*}
G_D(x,A)&=\E_x\left[\int_0^\infty \mathbbm{1}_{\{X^D_t\in A\}}\,dt\right]\\
&=\E_x\left[\int_0^\infty \mathbbm{1}_{\{X_t\in
A\}}\,dt\right]-\E_x\left[\int_{\tau_D}^\infty \mathbbm{1}_{\{X_t\in
A\}}\,dt\right]\\
&=G(x,A)-\E_x[G(X_{\tau_D},A); \tau_D<\infty]\,.
\end{align*}
Thus in the transient case the Green function of $X^D$ can be written as
\[
	G_D(x,y)=G(x,y)-\E_x[G(X_{\tau_D},y);\tau_D<\infty],\ x,y\in D,\
x\not=y\,.
\]

Since
$X$ is, in particular, an isotropic L\' evy process it follows from \cite{Sz} that 
\[
	\P_x(X_{\tau_{B_r(0)}}\in \partial B_r(0))=0,\ x\in B_r(0)\,.
\]
for any $r>0$ and $x\in B_r(0)$ .
This allows us to use the Ikeda-Watanabe formula (see   Theorem 1 in
\cite{IW}): 
\begin{equation}\label{eq:iw1}
	\P_x(X_{\tau_{B_r(0)}}\in
F)=\int_F\int_{B_r(0)}G_{B_r(0)}(x,y)j(|z-y|)\,dy\,dz\,,
\end{equation}
for $x\in B_r(0)$ and $F\subset \overline{B_r(0)}^c\,$.

Defining a function $K_{B_r(0)}\colon B_r(0)\times
\overline{B_r(0)}^c\rightarrow [0,\infty)$ by
\begin{equation}\label{eq:iw2}
	K_{B_r(0)}(x,z)=\int_{B_r(0)}G_{B_r(0)}(x,y)j(|z-y|)\,dy
\end{equation}
 the Ikeda-Watanabe formula (\ref{eq:iw1}) reads
\begin{equation}\label{eq:tmp_repr_intro}
\P_x(X_{\tau_{B_r(0)}}\in F)=\int_F K_{B_r(0)}(x,z)\,dz\,.
\end{equation}
The function $K_{B_r(0)}$ will be called the Poisson kernel for the ball
$B_r(0)\,$.

\subsection{Bernstein functions and subordinators}

A function $\phi\colon (0,\infty)\rightarrow (0,\infty)$ is said to be a
Bernstein function if $\phi\in C^\infty (0,\infty)$ and $(-1)^n\phi^{(n)}\leq 0$
for all $n\in \N$. Every Bernstein $\phi$ function has the following
representation:
\begin{equation}\label{eq:tmp_bernst_1}
	\phi(\lambda)=\gamma_1+\gamma_2 \lambda +\int_{(0,\infty)}(1-e^{-\lambda
t})\mu(dt),
\end{equation}
where $\gamma_1,\gamma_2\geq 0$ and $\mu$ is a measure on $(0,\infty)$
satisfying $\int_{(0,\infty)}(1\wedge t)\mu(dt)<\infty$. 

Using the elementary inequality $ye^{-y}\leq 1-e^{-y}$, $y>0$ we deduce from 
(\ref{eq:tmp_bernst_1}) that every Bernstein function $\phi$ satisfies:
\begin{equation}\label{eq:tmp-bernst-ineq}
\lambda\phi'(\lambda)\leq \phi(\lambda)\ \text{ for any }\ \lambda >0\,.
\end{equation}

There is a strong connection between subordinators and Bernstein functions. 
To be more precise,  $\phi$ is a Bernstein function such that $\phi(0+)=0$ (i.e.
$\gamma_1=0$ in
(\ref{eq:tmp_bernst_1})) if and only if it is the Laplace exponent of some
subordinator. If $\phi(0+)>0$, $\phi$ can be understood as the Laplace
exponent of a subordinator killed with rate $\phi(0+)$ (see Chapter 3 in
\cite{Be}).

A Bernstein function $\phi$ is a complete Bernstein function if the L\' evy
measure in (\ref{eq:tmp_bernst_1}) has a completely monotone density, i.e.
$\mu(dt)=\mu(t)\,dt$, where $\mu\colon (0,\infty)\rightarrow (0,\infty)$,
$\mu\in C^\infty (0,\infty)$ and $(-1)^n\mu^{(n)}\geq 0$ for any $n\in \N$ .

Let us mention some properties of complete Bernstein function (see \cite{SSV}).
The composition of two complete Bernstein function is a complete Bernstein
function and if $\phi$ is a complete Bernstein function, then
$\phi^\star(\lambda)=\frac{\lambda}{\phi(\lambda)}$ is also a complete Bernstein
function. 

Assume that $S$ is a subordinator with the Laplace exponent $\phi$ and infinite
L\'
evy measure. Then $\phi^\star$ is a Bernstein function if and only if the
potential measure $U$ has a decreasing density $u$ with respect to the Lebesgue
measure. Moreover, if $\nu$ denotes the L\' evy measure of the subordinator with
the Laplace exponent $\phi^\star$, then $u(t)=\nu(t,\infty)$ for any $t>0$ . 

\subsection{Regular variation}\label{subsec:reg_var}
A function $f\colon (0,\infty)\rightarrow (0,\infty)$ varies regularly 
(at infinity) with index $\rho\in \R$ if \[\lim\limits_{\lambda \to
+\infty}\frac{f(\lambda x)}{f(\lambda)}=x^\rho\ \text{ for every }\ x>0\,.\]
 If $\rho=0$, then we say that $f$ is slowly varying. Regular (slow) variation
at $0$ is defined similarly.

If $f$ varies regularly with index $\rho\in\R$, then there exists a slowly
varying function $\ell$ so that $f(\lambda)=\lambda^\rho \ell(\lambda)$\,.

Let $\ell$ be a slowly varying function such that $L(\lambda)=\int_0^\lambda
\frac{\ell(t)}{t}\,dt$ exists for all $\lambda>0$\,. Then $L$ is slowly varying,
\begin{equation}\label{eq:de_haan_limit}
 \lim_{\lambda\to+\infty }\frac{L(\lambda)}{\ell(\lambda)}=+\infty
\end{equation}

and 
\[
 \lim_{\lambda \to +\infty}\frac{L(\lambda x)-L(\lambda)}{\ell(\lambda)}=\log
x\ \text{ for every }\ x>0\,.
\]
(see Proposition 1.5.9 a and p. 127 in \cite{BGT}).

\section{Subordinators}\label{sec:asymp}

Let $S=(S_t\colon t\geq 0)$ be a subordinator with the Laplace
exponent $\phi$  satisfying the following conditions:
\begin{itemize}
 \item[{\bf (A-1)}]  there is $\alpha\in [0,2]$ such that  $\phi'$
varies regularly at infinity with index $\frac{\alpha}{2}-1$, i.e.
\[
  \lim_{\lambda \to+\infty }\frac{\phi'(\lambda
x)}{\phi'(\lambda)}=x^{\frac{\alpha}{2}-1}\ \text{ for every } \ x>0\,,
\]
\item[{\bf (A-2)}] the L\' evy measure is infite and has a decreasing
density $\mu$ ,
\item[{\bf (A-3)}] the potential measure has a decreasing density $u$\,.
\end{itemize}
When $\alpha=2$ we additionaly assume:
\begin{itemize}
 \item[{\bf (A-4)}] $\lambda \mapsto
\left(\frac{\lambda}{\phi(\lambda)}\right)'$ varies regularly at infinity with
index $-1$\,.
\end{itemize}

\begin{Rem}
\begin{itemize}
 \item[(a)] The most important assumption is (A-1) (and (A-4) when
$\alpha=2$). Other assumptions hold for a large class of subordinators
(e. g. when $\phi$ is a complete Bernstein function with infinite L\' evy
measure).
\item[(b)] By Karamata's theorem (see Theorem 1.5.11 in  \cite{BGT}) it follows
from  (A-1) that $\phi$ varies regularly at infinity with index
$\frac{\alpha}{2}$ .
\end{itemize}

\end{Rem}

\begin{Prop}\label{prop:sub-mu1}
	Let $\alpha\in[0,2)$ and let $S$ be a subordinator satisfying (A-1) and (A-2). Then 
	\[
		\mu(t)\asymp t^{-2}\phi'(t^{-1}),\ t\to 0+\,.
	\]	
\end{Prop}
\proof
Let $\varepsilon>0$. By a change of variable
\begin{align}
	\phi(\lambda+\varepsilon)-\phi(\lambda)&=\int_0^\infty (e^{-\lambda t}-e^{-(\lambda+\varepsilon)t})\mu(t)\,dt\nonumber \\
	&=\lambda^{-1}\int_0^\infty e^{-t}(1-e^{-\varepsilon \lambda^{-1}t})\mu(\lambda^{-1} t)\,dt\,.\label{eq:tmp_sub-1}
\end{align}
Since $\mu$ is decreasing, the following holds
\[
	\phi(\lambda+\varepsilon)-\phi(\lambda)\geq \lambda^{-1}\mu(\lambda^{-1})\int_0^1 e^{-t}(1-e^{-\varepsilon \lambda^{-1}t})\,dt\,.
\]
Now we can apply Fatou lemma to deduce
\begin{align*}
	\phi'(\lambda)&=\lim_{\varepsilon \to 0+}\frac{\phi(\lambda +\varepsilon)-\phi(\lambda)}{\varepsilon}\geq \lambda^{-2}\mu(\lambda^{-1})\int_0^1te^{-t}\,dt\\
	&=\lambda^{-2}\mu(\lambda^{-1} )(1-2e^{-1})\,.
\end{align*}
By setting $\lambda=t^{-1}$ we get the upper bound
\begin{equation}\label{eq:upper_bound_mu}
	\mu(t)\leq \frac{t^{-2}f'(t^{-1})}{1-2e^{-1}}\ \text{ for every }\ t>0\,.
\end{equation}

Now we prove the lower bound. Using  (\ref{eq:tmp_sub-1}), for any $r\in (0,1)$ and $\varepsilon>0$ we can write
\begin{equation}\label{eq:tmp-sun-1002}
\phi(\lambda+\varepsilon)-\phi(\lambda)=I_1+I_2
\end{equation}
 with
 \begin{align*}
 	I_1&=\lambda^{-1}\int_0^r e^{-t}(1-e^{-\varepsilon \lambda^{-1} t})\mu(\lambda^{-1}t)\,dt\\
	I_2&=\lambda^{-1}\int_r^\infty e^{-t}(1-e^{-\varepsilon \lambda^{-1} t})\mu(\lambda^{-1}t)\,dt\, .
 \end{align*}
Since $\mu$ is decreasing, the dominated convergence theorem yields
 \begin{align}
 	\limsup_{\varepsilon \to 0+}\frac{I_2}{\varepsilon}&\leq
\lambda^{-2}\mu(\lambda^{-1} r)\int_r^\infty te^{-t}\,dt\nonumber\\
	&=(r+1)e^{-r}\lambda^{-2}\mu(\lambda^{-1} r)\,.\label{eq:tmp-sun-1003}
 \end{align}
 
 To handle $I_1$  first we use the theorem of Potter (see Theorem 1.5.6
(iii) in \cite{BGT}) to conclude that there are constants $c_1>0$ and $\delta>0$
such that
 \begin{equation}\label{eq:tmp-sub-1001}
 	\frac{\phi'(\lambda t^{-1})}{\phi'(\lambda)}\leq c_1t^{\delta}\ \text{
for all }\ \lambda \geq 1\ \text{ and }\ t\leq 1\,.
 \end{equation}
Therefore, by (\ref{eq:upper_bound_mu}), (\ref{eq:tmp-sub-1001}) and the
dominated convergence theorem
\begin{align}
 \limsup_{\varepsilon \to 0+}\frac{I_1}{\varepsilon}&\leq \limsup_{\varepsilon
\to 0+}\frac{1}{1-2e^{-1}}\int_0^r e^{-t}\,\frac{1-e^{-\varepsilon
\lambda^{-1}t}}{\varepsilon\lambda^{-1} t}\,t^{-1}\phi'(\lambda
t^{-1})\,dt\nonumber\\
&\leq \phi'(\lambda)\,\frac{c_1}{1-2e^{-1}}\int_0^r
t^{\delta-1}e^{-t}\,dt\,.\label{eq:tmp-sun-1004}
\end{align}
Combining (\ref{eq:tmp-sun-1002}), (\ref{eq:tmp-sun-1003}) and
(\ref{eq:tmp-sun-1004}) we deduce
\[
 \phi'(\lambda)\leq \phi'(\lambda)\,\frac{c_1}{1-2e^{-1}}\int_0^r
t^{\delta-1}e^{-t}\,dt+(r+1)e^{-r}\lambda^{-2}\mu(\lambda^{-1} r)\,.
\]
Furthermore, by choosing $r\in (0,1)$ so that $\frac{c_1}{1-2e^{-1}}\int_0^{r}
t^{\delta-1}e^{-t}\,dt\leq \frac{1}{2}$ we get
\[
 \mu(\lambda^{-1}r)\geq \frac{e^r}{2(r+1)}\lambda^2\phi'(\lambda)\ \text{ for
every }\ \lambda \geq 1\,.
\]
By (A-1) we can find $t_0\in (0,r)$ so that
$\frac{\phi'(rt^{-1})}{\phi'(t^{-1})}\geq \frac{r^{\frac{\alpha}{2}-1}}{2}$ for
any $t\in (0,t_0)$. 
The lower bound nwo follows:
\[
 \mu(t)\geq \frac{r^{1+\frac{\alpha}{2}}e^r}{4(r+1)}t^{-2}\phi'(t^{-1})\ \text{
for every }\ t\in (0,t_0)\,.
\]
\qed

\begin{Rem}
The precise asymptotical behavior of $\mu$ when $\alpha\in (0,2)$  can be
obtained by Karamata's Tauberian theorem. This is not the case when
$\alpha=0$.\\
Under an additional
assumption
\[
 t\mapsto ta^{at}\mu(t)\ \text{ is
monotone on }\ (0,T) \ \text{ for some }\ a\geq 0\ \text{ and }\ T>0,
\]
it is possible to prove the following precise asymptotics in the case $\alpha=0$ :
\[
\mu(t)\sim t^{-2}\phi'(t^{-1}),\
t\to 0+\,.
\]
\end{Rem}

\begin{Prop}\label{prop:pot-dens}
 Let $\alpha\in[0,2)$ and let $S$ be a subordinator satisfying (A-1) and (A-3).
Then 
	\[	
u(t)\sim\frac{1}{\Gamma\left(1-\frac{
\alpha}{2}\right)}\frac{t^{-2}\phi'(t^{-1})}{\phi(t^{-1})^2},\ t\to 0+\,.
	\]	
\end{Prop}
\begin{Rem}
It can be proved that
\[	
u(t)\asymp\frac{1}{\Gamma\left(1-\frac{
\alpha}{2}\right)}\frac{t^{-2}\phi'(t^{-1})}{\phi(t^{-1})^2},\ t\to 0+\,.
	\]
similarly as in Proposition
\ref{prop:sub-mu1}. It is enough to note
\[
 \psi(\lambda+\varepsilon)-\psi(\lambda)=\int_0^\infty (e^{-\lambda
t}-e^{-(\lambda+\varepsilon)t}) u(t)\,dt
\]
with $\psi(\lambda)=-\frac{1}{\phi(\lambda)}$.

The main reason why we need precise asymptotics of $u$ is to be able to handle
the
case $\alpha=2$ by duality. 
\end{Rem}

\proof[Proof of Proposition \ref{prop:pot-dens}]
Let us first consider the case $\alpha=0$. In this case
$\ell(\lambda)=\lambda\phi'(\lambda)$ varies slowly (at infinity) and thus it
follows from Subsection \ref{subsec:reg_var} that
$\phi(\lambda)=\int_0^\lambda \frac{\ell(t)}{t}\,dt$ also varies slowly and 
\[
\lim_{\lambda\to +\infty} \frac{\phi(\lambda x)-\phi(\lambda)}{\lambda \phi'(\lambda)}=\log x\ \text{
for every }\ x>0\,.
\]

This and   (\ref{eq:prelim-tmp1}) imply
	\[
		\frac{\mathcal{L}U\left(\frac{1}{\lambda
x}\right)-\mathcal{L}U\left(\frac{1}{\lambda}\right)}{\frac{1}{\lambda}
\phi'\left(\frac{1}{\lambda}\right)\phi\left(\frac{1}{\lambda}\right)^2}=\frac{
\phi\left(\frac{1}{\lambda}\right)-\phi\left(\frac{1}{\lambda
x}\right)}{\frac{1}{\lambda}\phi'\left(\frac{1}{\lambda}\right)}\frac{
\phi\left(\frac{1}{\lambda}\right)}{\phi\left(\frac{1}{\lambda x}\right)}\to
\log x,\ \lambda\to 0+
	\]
	for any $x>0$. Now we can apply de Haan's Tauberian theorem (see
$0\,$--version of \cite[Theorem 3.9.1]{BGT}) to deduce
	\[
		\frac{U(\lambda
x)-U(\lambda)}{\frac{1}{\lambda}\phi'\left(\frac{1}{\lambda}
\right)\phi\left(\frac{1}{\lambda}\right)^2}\to \log x,\ \lambda\to 0+\,.
	\]
	
	If we apply de Haan's monotone density theorem (see  \cite[Theorem
3.6.8]{BGT}) we finally obtain
	\[
		u(t)\sim
\frac{t^{-2}\phi'\left(t^{-1}\right)}{t\phi\left(t^{-1}\right)^2},\ t\to 0+\,.
	\]

The case $\alpha\in (0,2)$ is already known. We give the proof for the sake of
completeness and adapt the result to the formula obtained in the case
$\alpha=0$ .

Since 
\[
	\mathcal{L}U(\lambda)=\frac{1}{\phi(\lambda)}
\]
varies regularly at infinity with index $-\frac{\alpha}{2}$,  Karamata's
Tauberian
theorem (see Theorem 1.7.1 in \cite{BGT})  implies
\[
	U([0,t])\sim
\frac{1}{\Gamma\left(1-\frac{\alpha}{2}\right)}\frac{1}{\phi\left(t^{-1}\right)}
, \
t\to 0+\,.
\]
Then by  Karamata's monotone density theorem (see Theorem 1.7.2 in \cite{BGT})
we deduce
\[
	u(t)\sim
\frac{\alpha}{2\Gamma\left(1-\frac{\alpha}{2}\right)}\frac{1}{
t\phi\left(t^{-1}\right)}\sim
\frac{1}{\Gamma\left(1-\frac{\alpha}{2}\right)}\frac{t^{-2}\phi'(t^{-1})}{
\phi\left(t^{-1}\right)^2},\ t\to 0+\,.
\]
\qed

Now we consider the case $\alpha=2$. 


 \begin{Prop}\label{prop:alpha2_m}  Let $\alpha=2$ and let $S$ be a
subordinator satisfying (A-1), (A-2) and (A-4).
Then 
 \[
 	\mu(t)\sim t^{-2}\left(t\phi(t^{-1})-\phi'(t^{-1})\right)\,,\ t\to
0+\,.
 \]	
 \end{Prop}
 
\proof	
	
	Since potential density exists by (A-3) we see that $\phi$ it follows
that  $\phi^\star(\lambda)=\frac{\lambda}{\phi(\lambda)}$ defines
the Laplace exponent of a (possibly killed) subordinator, which we denote by
$T$ (see Section \ref{sec:prelim}).
	
 Note that the subordinator $T$ corresponds to the case of $\alpha=0$.  
If we denote potential density of $T$ by $v$, then
 \[
 	v(t)=\mu(t,\infty),\ t>0\,.
 \]
Proposition \ref{prop:pot-dens} yields
 \begin{equation}\label{eq:tmp20111}
 	\int_t^\infty \mu(s)\,ds\sim
\frac{t^{-2}(\phi^\star)'(t^{-1})}{\phi^\star(t^{-1})^2},\ t\to 0+\,.
 \end{equation}
By assumption (A-4) we know that $t\mapsto (\phi^\star)'(t^{-1})$ varies
regularly at
$0$ with index $1$ and thus
$t\mapsto \frac{t^{-2}(\phi^\star)'(t^{-1})}{\phi^\star(t^{-1})^2}$ varies
regularly at
$0$ with index $-1$. 

Now we change variable in the integral on the left-hand side in
(\ref{eq:tmp20111}) and conclude
\[
	\int_0^{t^{-1}}\mu(s^{-1})\,\frac{ds}{s^2}\sim
\frac{t^{-2}(\phi^\star)'(t^{-1})}{(\phi^\star(t^{-1}))^2},\ t\to 0+\,.
\]
This gives ($r=t^{-1}$):
\[
	\int_0^{r}\mu(s^{-1})\,\frac{ds}{s^2}\sim
\frac{r^2(\phi^\star)'(r)}{\phi^\star(r)^2},\ r\to \infty\,.
\]
Note that the right-hand side is now regularly varying at infinity with index
$1$ and thus by Karamata's monotone density theorem (see Theorem 1.7.2 in
\cite{BGT}) we deduce
\[
	\frac{\mu(r^{-1})}{r^2}\sim
\frac{r(\phi^\star)'(r)}{\phi^\star(r)^2}\,.
r\to\infty\,,
\]
Going back ($t=r^{-1}$) we conclude
\[
 \mu(t)\sim \frac{t^{-3}(\phi^\star)'(t^{-1})}{\phi^\star(t^{-1})^2},\ t\to
0+\,.
\]

\qed

  \begin{Prop}\label{thm:alpha2_u}  Let $\alpha=2$ and let $S$ be a
subordinator satisfying (A-1) and  (A-3). Then the following is true
 \[
 	u(t)\sim \frac{1}{\phi'(t^{-1})}\sim  \frac{1}{t\phi(t^{-1})},\ t\to
0+\,.
 \]	
 \end{Prop}
\proof
	By (\ref{eq:prelim-tmp1}) we get
	\[
		\mathcal{L}U(\lambda)=\frac{1}{\phi(\lambda)}\sim
\frac{1}{\lambda\phi'(\lambda)},\ \lambda\to +\infty
	\]
	and thus by  Karamata's Tauberian theorem (see Theorem 1.7.1 in
\cite{BGT}) it follows that 
	\[
		U([0,t])\sim \frac{1}{\Gamma(2)}\frac{t}{\phi'(t^{-1})},\ t\to
0+
	\]	
	since $\lambda\mapsto\lambda\phi'(\lambda)$ varies regualrly at infinity
with index $1$.
	By applying Karamata's monotone density (see Theorem 1.7.2 in
\cite{BGT}) theorem we deduce
	\[
		u(t)\sim \frac{1}{\phi'(t^{-1})},\ t\to 0+\,.
	\]
\qed

\section{L\' evy density and Green function}\label{sect:sbm}

Let $S$ be a subordinator as in Section \ref{sec:asymp} and let $X$ be the
corresponding subordinate Brownian motion in $\R^d$ with $d\geq 3$. Our aim is
to establish asymptotical behavior of the L\' evy density and Green function of
$X$ .

Recall that the L\' evy density of $X$ is of the form $j(|x|)$, where $j$ is
given by (\ref{eq:prelim-jot}).

\begin{Thm}\label{thm:jump_asymp}
 Assume that $S$ satisfies (A-1) with some $\alpha\in [0,2]$ and (A-2). If
$\alpha\in [0,2)$, then
\[
 j(r)\asymp r^{-d-2}\phi'(r^{-2}),\ r\to 0+\,.
\]
If $\alpha=2$ and (A-4) holds, then
\[
 j(r)\asymp r^{-d-2}\left(r^2\phi(r^{-2})-\phi'(r^{-2})\right),\ r\to 0+\,.
\]
\end{Thm}
\proof
  This result follows directly from Proposition \ref{prop:sub-mu1} and
Proposition \ref{prop:alpha2_m}  together with Lemma \ref{lem:asympt_lema}, 
where $a=\frac{1}{4}$, $b=1+\frac{\alpha}{2}$. 

For $\alpha\in [0,2)$ the slowly
varying function is given by 
$\ell(t)=t^{\frac{\alpha}{2}-1}\phi'(t^{-1})$. When $\alpha=2$ we take

\[
 \ell(t)=t\phi(t^{-1})-\phi'(t^{-1})=\frac{t\left(\phi^\star(t^{-1}
)\right)' } {\phi^\star(t^ { -1})^2},
\]
which varies slowly, 
since $\phi^\star(\lambda)=\frac{\lambda}{\phi(\lambda)}$ is slowly varying by
(A-1) and Karamata's Theorem (see Theorem 1.5.11 in \cite{BGT})  and has a
derivative that varies regularly with index $1$ by (A-4).
\qed

The Green function of $X$ is of the form $G(x,y)=g(|y-x|)$, where $g$ is given
by (\ref{eq:prelim-green}).

\begin{Thm}\label{thm:green_asymp}
 Assume that $S$ satisfies (A-1) with some $\alpha\in [0,2]$ and (A-3). Let
$d\geq 3$. If
$\alpha\in [0,2)$, then
\[
 g(r)\asymp r^{-d-2}\frac{\phi'(r^{-2})}{\phi(r^{-2})^2},\ r\to 0+\,.
\]
If $\alpha=2$,  then
\[
 g(r)\asymp r^{-d+2}\frac{1}{\phi'(r^{-2})}\asymp r^{-d}\frac{1}{\phi(r^{-2})},\
r\to 0+\,.
\]
\end{Thm}
\proof
  We use Lemma \ref{lem:asympt_lema} with $a=\frac{1}{4}$, 
$b=1-\frac{\alpha}{2}$, Proposition \ref{prop:pot-dens} and Proposition
\ref{thm:alpha2_u}.

When $\alpha\in [0,2)$ we define
$\ell(t)=\frac{t^{-1}\phi'(t^{-1})}{\phi(t^{-1})}$ which varies slowly at $0$ .

In the case $\alpha=2$, we let $\ell(t)=\frac{1}{\phi'(t^{-1})}$ or
$\ell(t)=\frac{1}{t\phi(t^{-1})}$  which both vary slowly at $0$.
\qed

Using  asymptotical results from this section we can now prove the proposition
that gives a counterexample for the estimate of the Krylov and Safonov.

\proof[Proof of Proposition \ref{prop:count-ks}]
  By (\ref{eq:iw1}),
\begin{align*}
 \P_x(X_{\tau_{B_{\frac{r}{4}}(0)}}\in B_r(0)\setminus
B_{\frac{r}{4}}(0))&=\int_{B_r(0)\setminus
\overline{B_{\frac{r}{4}}(0)}}\int_{B_{\frac{r}{4}}(0)}G_{B_{\frac{r}{4}}}(x,
y)j(|z-y|)\, dy\\&=I_1+I_2\,.
\end{align*}

Using Theorems \ref{thm:jump_asymp} and  \ref{thm:green_asymp} it follows that
\begin{align*}
 I_1&=\int_{B_r(0)\setminus
\overline{B_{\frac{r}{4}}(0)}}\int_{B_{\frac{3r}{16}}(0)}G_{B_{\frac{r}{4}}}(x,
y)j(|z-y|)\,dy\,dz\\
&
\leq j(\tfrac{r}{16})|B_r(0)\setminus
B_{\frac{r}{4}}(0)|\int_{B_r(0)}g(|y|)\,dy\,dz\\
&\leq c_1 \frac{r^{-2}\phi'(r^{-2})}{\phi(r^{-2})}\,.
\end{align*}
On the other hand,
\begin{align}
 I_2&=\int_{B_r(0)\setminus
\overline{B_{\frac{r}{4}}(0)}}\int_{B_{\frac{r}{4}}(0)\setminus
B_{\frac{3r}{16}}
(0) } G_ { B_ { \frac { r } { 2 } } } (x ,
y)j(|z-y|)\,dy\,dz\nonumber\\
&\leq g(\tfrac{r}{16})\int_{B_r(0)\setminus
\overline{B_{\frac{r}{4}}(0)}}\int_{B_{\frac{r}{4}}(z)}j(|y|)\,dy\,dz\,.\label{eq:count-tmp11}
\end{align}
To estimate the inner integral, note that $B_{\frac{r}{4}}(z)\subset
B_1(0)\setminus B_{|z|-\frac{r}{4}}(0)$ for any $z\in B_r(0)\setminus
\overline{B_{\frac{r}{4}}(0)}$ and so, by Theorem \ref{thm:jump_asymp}, 
\begin{equation}\label{eq:jump_gr-tmp10}
 \int_{B_{\frac{r}{4}}(z)}j(|y|)\,dy\leq
c_2\int_{|z|-\frac{r}{4}}^1s^{-3}\phi'(s^{-2})\,ds\leq c_3
\phi((|z|-\tfrac{r}{4})^{-2})\,.
\end{equation}

Thus, by Theorem \ref{thm:green_asymp} and (\ref{eq:count-tmp11})
\begin{align*}
 I_2&\leq  c_4
r^{-d-2}\frac{\phi'(r^{-2})}{\phi(r^{-2})^2}\int_{\frac{r}{4}}
^r\phi((t-\tfrac{r}{4})^{-2})t^{d-1}\,dt\\
&\leq c_4
r^{-3}\frac{\phi'(r^{-2})}{\phi(r^{-2})^2}\int_0^{\frac{r}{4}}\phi(s^{-2})\,ds\\
&\leq c_5 \frac{r^{-2}\phi'(r^{-2})}{\phi(r^{-2})}.
\end{align*}
In the last equality we have used Karamata's theorem (see Theorem 1.5.11 in
\cite{BGT}) and the fact that $\alpha\in [0,1)$ .
\qed

\pagebreak

\section{Difference estimates}\label{sect:diff}

Let $X$ be the stochastic process in $\R^d$ as in Section \ref{sect:sbm} and
assume that $d\geq 3$. In particular $X$ is transient.

In this section we prove the difference estimates of the Green function
and the Poisson kernel.

Although we are slightly abusing notation, we set $G(x):=G(0,x)=g(|x|)$. 

\begin{Prop}\label{prop:green_diff_1}
	There is a constant $c>0$ such that for every $r\in (0,1)$
	\[
		|G(x)-G(y)|\leq cg(r)\left(1\wedge \tfrac{|x-y|}{r}\right)\
\text{ for  all } \ x,y\not\in
B_r(0)\,.
	\]
\end{Prop}
\proof
	Assume first that $|x-y|<\frac{r}{2}$. By the mean value theorem it
follows that
for any $t>0$ there exists $\vartheta=\vartheta(x,y,t)\in [0,1]$ such that 
	\begin{align*}
		\left|e^{-\frac{|x|^2}{4t}}-e^{-\frac{|y|^2}{4t}}\right|&\leq
\tfrac{|x+\vartheta(y-x)|}{2t}e^{-\frac{|x+\vartheta(y-x)|^2}{4t}}|x-y|\\
		&\leq
2\tfrac{|x-y|}{\sqrt{t}}e^{-\frac{|x+\vartheta(y-x)|^2}{8t}}\,,
	\end{align*}
	where in the last line the following elementary inequality was used
	\[
		se^{-s^2}<2 e^{-\frac{s^2}{2}},\ s>0\,.
	\]
	Then
	$
		|x+\vartheta(y-x)|\geq |x|-\vartheta|y-x|\geq \frac{r}{2}
	$
	implies
	\begin{equation}\label{eq:green_diff_1}
		\left|e^{-\frac{|x|^2}{4t}}-e^{-\frac{|y|^2}{4t}}\right|\leq
2\tfrac{|x-y|}{\sqrt{t}}e^{-\frac{r^2}{32t}}\,.
	\end{equation}
	By (\ref{eq:green_diff_1})
	\begin{align*}
		|G(x)-G(y)|&\leq (4\pi)^{-d/2}\int_0^\infty
t^{-d/2}\left|e^{-\frac{|x|^2}{4t}}-e^{-\frac{|y|^2}{4t}}\right|u(t)\,dt\\
		&\leq 2(4\pi)^{-d/2}|x-y|\int_0^\infty
t^{-d/2-1/2}e^{-\frac{r^2}{32t}}u(t)\,dt\,.
	\end{align*}
	Since $u$ is non-increasing and varies regularly at $0$ with index
$\frac{\alpha}{2}-1$, by  Lemma \ref{lem:asympt_lema} we see that there is a
constant
$c_1>0$ so that 
	\[
		\int_0^\infty t^{-d/2-1/2}e^{-\frac{r^2}{32t}}u(t)\,dt\leq c_1
r^{-d+1}u(r^2)\ \text{ for every }\ r\in (0,1)\,.
	\]
	Theorem \ref{thm:green_asymp} yields
	\[
		|G(x)-G(y)|\leq c_2\,g(r)\tfrac{|x-y|}{r}\,.
	\]
	
	When $|x-y|\geq \frac{r}{2}$, 
	\begin{align*}
		|G(x)-G(y)|&\leq G(x)+G(y)\leq 2g(r)
	\end{align*}
	since $|x|,|y|\geq r\,$.
\qed

\begin{Prop}\label{prop:green_diff_2}
	There is a constant $c>0$ such that for all $R\in (0,1)$, $r\in (0,
\frac{R}{2}]$,
$y\in B_R(0)$ and $x_1,x_2\in B_{\frac{R}{2}}(0)\setminus B_r(y)$
	\[
		|G_{B_R(0)}(x_1,y)-G_{B_R(0)}(x_2,y)|\leq c g(r)\left(1\wedge
\tfrac{|x_1-x_2|}{r}\right)\,.
	\]
\end{Prop}
\proof
	By symmetry of the Green function, 
	\begin{align*}	
G_{B_R(0)}(x_i,y)&=G_{B_R(0)}(y,x_i)=G(x_i-y)-\E_y[G(X_{\tau_{B_R(0)}}-x_i)]
\\&=G(x_i-y)-\E_y[G(X_{\tau_{B_R(0)}}-x_i)]\,,
	\end{align*}
	for $i\in \{1,2\}\,$.
	 Now the result follows from Proposition \ref{prop:green_diff_1}.
\qed


\begin{Prop}\label{prop:poisson_k_est}
	There is a constant $c>0$ such that for any $r\in (0,1)$ and  $x,y\in
B_{\frac{r}{8}}(0)$:
	\begin{itemize}
	\item[(i)] if $z\in B_{2r}(0)\setminus B_r(0)$, then
	\[
		\left|K_{B_r(0)}(x,z)-K_{B_r(0)}(y,z)\right|\leq
c|z|^{-d}\frac{\phi\left((|z|-r)^{-2}\right)}{\phi\left(|x-y|^{-2}\right)}
\,;
	\]
	\item[(ii)] if $z\not\in B_{2r}(0)$, then
	\[
		\left|K_{B_r(0)}(x,z)-K_{B_r(0)}(y,z)\right|\leq
c\frac{j\left(\frac{|z|}{2}\right)}{\phi\left(|x-y|^{-2}\right)}
\,.	\]
	\end{itemize}
\end{Prop}
\proof
In the estimate
\begin{align*}
	\left|K_{B_r(0)}(x,z)-K_{B_r(0)}(y,z)\right|&\leq
\int_{B_r(0)}\left|G_{B_r(0)}(x,v)-G_{B_r(0)}(y,v)\right|j(|z-v|)\,dv
\end{align*}
we split the integral into three parts:
\begin{align*}
	I_1=&\int_{B_{2|x-y|}(x)}\left|G_{B_r(0)}(x,v)-G_{B_r(0)}(y,
v)\right|j(|z-v|)\,dv\\
	I_2=&\int_{B_{\frac{r}{4}}(x)\setminus
B_{2|x-y|}(x)}\left|G_{B_r(0)}(x,v)-G_{B_r(0)}(y,v)\right|j(|z-v|)\,dv\\
	I_3=&\int_{B_{r}(0)\setminus
B_{\frac{r}{4}}(x)}\left|G_{B_r(0)}(x,v)-G_{B_r(0)}(y,v)\right|j(|z-v|)\,dv\,.
\end{align*}
For the first part we obtain
\begin{align}
	I_1&\leq \int_{B_{2|x-y|}(x)}G_{B_r(0)}(x,v)j(|z-v|)\,dv+
\int_{B_{3|x-y|}(y)}G_{B_r(0)}(y,v)j(|z-v|)\,dv\nonumber\\
	&\leq 2j\left(\tfrac{|z|}{2}\right)\int_{B_{3|x-y|}(0)}G(v)\,dv\leq
c_1\tfrac{j\left(\frac{|z|}{2}\right)}{\phi\left(|x-y|^{-2}\right)}\,,
\label{eq:int_part1}
\end{align}
for any $z\not\in B_r(0)$. We have used Theorem \ref{thm:green_asymp}
to get the last inequality in (\ref{eq:int_part1}).

In order to estimate $I_2$ we split the integral in the following way. We let
$N=\left\lfloor \frac{\log\frac{r}{4|x-y|}}{\log 2}\right\rfloor$  and write
\[
	I_2\leq \sum_{n=1}^N\int_{B_{2^{n+1}|x-y|}(x)\setminus
B_{2^{n}|x-y|}(x)}\left|G_{B_r(0)}(x,v)-G_{B_r(0)}(y,v)\right|j(|z-v|)\,dv\,.
\]
Now, for each $n\in\{1,\ldots,N\}$ we can apply Proposition
\ref{prop:green_diff_2} (with the corresponding radii $(2^{n}-1)|x-y|$ and $r$) to get
\begin{align*}
	\int_{B_{2^{n+1}|x-y|}(x)\setminus
B_{2^{n}|x-y|}(x)}&\left|G_{B_r(0)}(x,v)-G_{B_r(0)}(y,v)\right|j(|z-v|)\,dv\\
	&\leq
c_3\frac{g\left((2^n-1)|x-y|\right)}{2^n-1}\int_{B_{2^{n+1}|x-y|}(x)}j(|z-v|)\,
dv\,.
\end{align*}
By Theorem \ref{thm:green_asymp}
\[
\frac{g\left((2^n-1)|x-y|\right)}{g(|x-y|)}\leq c_4
\frac{\eta((2^n-1)|x-y|)}{\eta(|x-y|)}\ \text{ for all }\ n\in \{1,2,\ldots,N\},
\]
with $\eta(r)=r^{-d-2}\frac{\phi'(r^{-2})}{\phi(r^{-2})^2}$ .

Noting
that 
$\eta$ varies regularly at zero
with index $\alpha-d<0$, the uniform convergence theorem for regularly
varying
functions (see Theorem 1.5.2 in  \cite{BGT}) gives
\[
	\frac{\eta\left((2^n-1)|x-y|\right)}{\eta(|x-y|)}\leq c_5
(2^n-1)^{\alpha-d}\
\ \textrm{ for all } \ n\in\N\ \text{ and } |x-y|\leq \tfrac{1}{2}\, .
\]

By Theorem \ref{thm:green_asymp} and (\ref{eq:tmp-bernst-ineq})
$
	g(|x-y|)\leq
\frac{c_5}{\phi(|x-y|^{-2})}
$
 and so
\begin{align}
	I_2&\leq c_6
\sum_{n=1}^N(2^n-1)^{\alpha-d-1}g(|x-y|)(2^{n+1}|x-y|)^dj\left(\tfrac{|z|}{
2}\right)\nonumber \\
	&\leq c_7
\frac{j\left(\frac{|z|}{2}\right)}{\phi\left(|x-y|^{-2}\right)}\sum_{n=1}^N
2^{(\alpha-1)n}\nonumber\\
&\leq \frac{c_7}{1-2^{\alpha-1}}
\frac{j\left(\frac{|z|}{2}\right)}{\phi\left(|x-y|^{-2}\right)} \ \text{ for
every }\ z\not\in B_r(0)\,.\nonumber
\end{align}


It remains to estimate $I_3$. Applying Theorem \ref{prop:green_diff_2} we
get
\begin{align}
	I_3&\leq c_{8} g(r)\frac{|x-y|}{r}\int_{B_r(z)}j(|v|)\,dv\nonumber\\
	    &\leq c_{9}
\frac{|x-y|\phi\left(|x-y|^{-2}\right)}{r\phi\left(r^{-2}\right)}\frac{r^{-d}}{
\phi\left(|x-y|^{-2}\right)}\int_{B_r(z)}j(|v|)\,dv\nonumber\\
	    &\leq c_{10}
\frac{r^{-d}}{\phi\left(|x-y|^{-2}\right)}\int_{B_r(z)}j(|v|)\,dv\,.\label{eq:tmp_102c}
\end{align}
In the last inequality we have used the theorem of Potter (cf.
\cite[Theorem 1.5.6 (iii)]{BGT}) to conclude that for
$\delta<1-\alpha$ there is a constant $A_\delta>0$ such that 
\[
	\frac{|x-y|\phi(|x-y|^{-2})}{r\phi(r^{-2})}\leq
A_\delta\left(\frac{|x-y|}{r}\right)^{1-\alpha-\delta}\leq
A_\delta,
\]
since $r\mapsto r\phi(r^{-2})$ varies regularly at zero with index $1-\alpha$ .

Since 
\[
	j(|v|)\geq j\left(\tfrac{|z|}{2}\right)\ \textrm{ for all }\ v\in B_r(z)\ \text{ and }\ z\in B_{2r}(0)^c
\]
it follows from (\ref{eq:tmp_102c}) that
\begin{align*}
	I_3&\leq c_{11}
\frac{j\left(\frac{|z|}{2}\right)}{\phi\left(|x-y|^{-2}\right)}\,,
\end{align*}
On the other hand, for $z\in B_{2r}(0)\setminus B_r(0)$ we deduce from 
\[
	B_r(z)\subset B_3(0)\setminus B_{|z|-r}(0)
\]
(similarly as in (\ref{eq:jump_gr-tmp10})) that
\begin{align*}
\int_{B_r(z)}j(|v|)\,dv&\leq c_{12} \phi\left((|z|-r)^{-2}\right)\,. 
\end{align*}
By (\ref{eq:tmp_102c}) 
\[
	I_3\leq
c_{13}|z|^{-d}\frac{\phi\left((|z|-r)^{-2}\right)}{\phi\left(|x-y|^{-2}
\right)}\ \text{ for all } z\in B_{2r}(0)\setminus B_r(0)\,.
\]
\qed
\section{Regularity of harmonic functions}\label{sect:main}


Recall that (\ref{eq:tmp_repr_intro}) gives the representation for any bounded
function $f\colon\R^d\rightarrow \R$ that is harmonic in $B_{2r}(x_0)$:
\begin{equation}\label{eq:harm_f_repr}
	f(x)=\E_x\left[f\left(X_{\tau_{B_r(x_0)}}\right)\right]=\int_{\overline{
B_r(x_0)}^c}K_{B_r(x_0)}(x,z)f(z)\,dz,\ \ x\in B_r(x_0)\,.
\end{equation}

\proof[Proof of Theorem \ref{tm:main}]
By  (\ref{eq:harm_f_repr})
\begin{equation}\label{eq:harm_tmp1}
	|f(x)-f(y)|\leq \|f\|_\infty
\int_{\overline{B_{2r}(0)}^c}\left|K_{B_{2r}(0)}(x,z)-K_{B_{2r}(0)}(y,
z)\right|\,dz\,.
\end{equation}
It remains to estimate the integral in (\ref{eq:harm_tmp1}), which we split in
the following way
\begin{align*}
	I_1&=\int_{B_{4r}(0)\setminus
\overline{B_{2r}(0)}}\left|K_{B_{2r}(0)}(x,z)-K_{B_{2r}(0)}(y,z)\right|\,dz\\
	I_2&=\int_{B_1(0)\setminus
B_{4r}(0)}\left|K_{B_{2r}(0)}(x,z)-K_{B_{2r}(0)}(y,z)\right|\,dz\\
	I_3&=\int_{B_1(0)^c}\left|K_{B_{2r}(0)}(x,z)-K_{B_{2r}(0)}(y,z)\right|\,
dz
\end{align*}
In order to estimate $I_1$ we use Proposition \ref{prop:poisson_k_est} (i). More
precisely,
\begin{align*}
	I_1&\leq \frac{c_1}{\phi\left(|x-y|^{-2}\right)}\int_{B_{4r}(0)\setminus
\overline{B_{2r}(0)}}|z|^{-d}\phi\left((|z|-2r)^{-2}\right)\,dz\\
	     &=
\frac{c_2}{\phi\left(|x-y|^{-2}\right)}\int_{2r}^{4r}t^{-1}\phi\left((t-2r)^{-2}
\right)\,dt\\&\leq \frac{c_2}{\phi\left(|x-y|^{-2}\right)}
(2r)^{-1}\int_0^{2r}\phi\left(t^{-2}\right)\,dt\\
	     &\leq \frac{c_3}{{\phi\left(|x-y|^{-2}\right)}} \phi(r^{-2})\,,
\end{align*}
where in the last inequality we have used Karamata's theorem (see the
$0$-version of  Theorem 1.5.11 in  \cite{BGT}).

We estimate $I_2$ and $I_3$ with the help of  Proposition
\ref{prop:poisson_k_est} (ii). Since the L\' evy measure is finite away from the
origin,
\[
	I_3\leq
\frac{c_4}{\phi\left(|x-y|^{-2}\right)}\int_{B_1(0)^c}j\left(\tfrac{|z|}{2}
\right)\,dz\leq\frac{c_5}{\phi\left(|x-y|^{-2}\right)}\,.
\]
Also, 
\begin{align*}
	I_2&\leq  \frac{c_6}{\phi\left(|x-y|^{-2}\right)}\int_{B_1(0)\setminus
B_{4r}(0)}j\left(\tfrac{|z|}{2}\right)\,dz
	\leq \frac{c_7\phi\left(r^{-2}\right)}{\phi\left(|x-y|^{-2}\right)}\,,
\end{align*}
where in the last inequality we have used Theorem \ref{thm:jump_asymp}.
\qed

\section{Examples}\label{sec:examples}

In this section is to illustrate our results by some
examples.

\subsection{(Iterated) Geometric stable processes}\label{subsec:geom_stable}

This class of examples belongs to the  case of $\alpha=0$.

Let $\beta\in (0,2]$. We define a family of functions $\{\phi_n\colon
(0,\infty)\rightarrow (0,\infty)\colon n\in \N\}$ recursively by
\begin{align*}
 \phi_1(\lambda)&=\log(1+\lambda^{\beta/2}),\ \lambda>0\\
\phi_{n+1}&=\phi_1\circ\phi_n,\ n\in \N\,.
\end{align*}

The function $\phi_1$ is a complete Bernsetin function. Since complete
Bernstein functions are closed under operation of composition, $\phi_n$ belongs
to this class for every $n\in \N$.

Let $S^n$ be a subordinator with the Laplace exponent $\phi_n$. $S^1$ is known
as the geometric $\frac{\beta}{2}$-stable subordinator. We call $S^n$
the iterated geometric $\frac{\beta}{2}$-stable subordinator. The corresponding
subordinate Brownian motions $X^n$ will be called (iterated) geometric
$\beta$-stable processes. 

As already remarked in \cite{SSV2}, these processes show quite different
behavior compared to the one of stable processes. Our contribution to this
class
of examples is that now we can obtain behavior of the L\' evy density as a
special case of Theorem \ref{thm:jump_asymp}(even for iterated geometric stable
processes).

The L\' evy density of $X^n$ is comparable to
\[
 \frac{1}{|x|^d}\cdot \prod_{k=1}^{n-1}\frac{1}{\log_{k}(|x|^{-1})}\ \text{ as
}\ |x|\to 0+\,,
\]
which is almost integrable. We can say that (intially) this
process jumps slower than any stable processes.

This can be also seen from the behavior of the Green function:
\[
 G(x,y)\asymp
\frac{1}{|x-y|^d\,\log_n^2(|x-y|^{-1})}\cdot \prod_{k=1}^{n-1}\frac{1}{
\log_{k}(|x-y|^{-1})}\ \text{ as
}\ |x-y|\to 0+\,\,.
\]
 As a consequence, $\E_0\tau_{B_r(0)}\asymp \frac{1}{\log_n (r^{-1})}$ as $r\to
0+$. Therefore $X^n$ needs (on average) more time to exit ball $B_r(0)$ than
any stable process or Brownian motion. 

Theorem \ref{tm:main} implies the  following a-priori local regularity 
estimates of harmonic functions:
\[
 |f(x)-f(y)|\leq c\|f\|_\infty \log_n(r^{-1})\,\frac{1}{\log_n(|x-y|^{-1})}\
\text{ for all }\ x,y\in B_{\frac{r}{2}}(0)
\]
and any bounded function $f$ which is harmonic in $B_r(0)$ . 

This tells us that the modulus of continuity is bounded by a logarithmic term.
It is still an open problem whether these harmonic functions satify a-priori
local H\" older continuity estimates.

\subsection{Conjugates of  (iterated) geometric stable
processes}\label{subsec:conj_geom_stable}

This class of examples corresponds to the case $\alpha=2$.

Let $\psi_n(\lambda)=\frac{\lambda}{\phi_n(\lambda)}$, where $\phi_n$ are as in
Subsection \ref{subsec:geom_stable}. 

Since $\phi_n$ are complete Bernstein
functions, $\psi_n$ are also complete Bernstein functions. Therefore, there
exist (killed) subordinators $T^n$ with the Laplace exponent $\psi_n$. Killing
will not affect the behavior of the L\' evy and potential densities of $T^n$
near zero. 

In this case the L\' evy density of the corresponding subordinate
Brownian motion $Y^n$ behaves near the origin as
\[
 \frac{1}{|x|^{d+2}\, \log_n^2(|x|^{-1})}\cdot
\prod_{k=1}^{n-1}\frac{1}{\log_k(|x|^{-1})} \ \text{ as }\ |x|\to 0+\,.
\]

Note that the integrability conditions of the L\' evy measure are barely
satisfied in this case. 

Comparing this behavior to the behavior of the small jumps of the
$\alpha$-stable process, we see that small jumps of $Y^n$ are more intensive.

Another interesting feature of this process is the following behavior of the
Green function:
\[
 G(x,y)\asymp |x-y|^{2-d}\log_n(|x-y|^{-1})\ \text{ as }\ |x-y|\to 0+\,.
\]
In this sense the process $Y^n$ is 'between' stable processes and Brownian
motion, since their Green functions are given by 
\[
  G^{(\alpha)}=c_\alpha|x-y|^{\alpha-d}\ \text{ and }\ 
G^{(2)}=c_\alpha|x-y|^{2-d}\,.
\]

\appendix
\section{Asymptotical properties}

In the appendix we prove a technical lemma which is used throughout the paper. 

\begin{Lem}\label{lem:asympt_lema}
	Let  $w\colon (0,\infty)\rightarrow (0,\infty)$ be a decreasing function
satisfying
	\[
		w(t)\asymp t^{-b}\ell(t),\ t\to 0+\,,
	\]
	for a function $\ell\colon (0,\infty)\rightarrow
(0,\infty)$ that varies slowly at $0$ and $b\geq 0\,$. 
	
	If $p>1$ and $a>0$, then 
	\[
		I(r)=\int_0^\infty t^{-p} e^{-\frac{ar}{t}}w(t)\,dt,\ r>0\,,
	\]
	satisfies
	\[
		I(r)\asymp a^{-p-b+1}r^{-p+1}w(r),\ r\to 0+\,.
	\]
\end{Lem}
\proof
Change variables yields
	\begin{equation}\label{eq:app-tmp1}	
I(r)=(ar)^{-p+1}\int_0^\infty
e^{-t}t^{p-2}w\left(\frac{ar}{t}\right)\,dt
	\end{equation}
	By assumptions, there are constants $c_1,c_2>0$ and
$r_0>0$ such that
\begin{equation}\label{eq:app-tmp2}
  c_1a^{-b}\leq \frac{w(ar)}{w(r)}\leq c_2 a^{-b} \ \text{ for every }\ r\in
(0,r_0)\,.
\end{equation}

	Let us first prove the upper bound. Using the fact that $w$ is
decreasing, by (\ref{eq:app-tmp1}) and (\ref{eq:app-tmp2}) we get
\begin{align*}
  I(r)&\leq (ar)^{-p+1}\int_0^1
e^{-t}t^{p-2}w\left(ar\right)\,dt+(ar)^{-p+1}\int_1^\infty
e^{-t}t^{p-2}w\left(\frac{ar}{t}\right)\,dt\\
&\leq c_2 (ar)^{-p+1}a^{-b}w(r)\left[\int_0^1 e^{-t}t^{p-2}\,dt+\int_1^\infty
e^{-t}t^{p+b-2}\,dt\right]\\
&\leq c_2'a^{-p-b+1}r^{-p+1}w(r)
\end{align*}
for every $r\in (0,r_0)$ .

The lower bound follows similarly:
\begin{align*}
  I(r)&\geq (ar)^{-p+1}\int_1^\infty
e^{-t}t^{p-2}w\left(ar\right)\,dt
\geq c_1(ar)^{-p+1}a^{-b}w(r)\int_1^\infty e^{-t}t^{p-2}\,dt\\
&= c_1'a^{-p-b+1}r^{-p+1}w(r)
\end{align*}
for every $r\in (0,r_0)$ .
\qed

\providecommand{\bysame}{\leavevmode\hbox to3em{\hrulefill}\thinspace}
\providecommand{\MR}{\relax\ifhmode\unskip\space\fi MR }
\providecommand{\MRhref}[2]{%
  \href{http://www.ams.org/mathscinet-getitem?mr=#1}{#2}
}
\providecommand{\href}[2]{#2}

\end{document}